\documentclass[11pt]{amsart}

\usepackage{amsmath, amssymb}
\usepackage{amsfonts}
\usepackage{amsthm}
\usepackage{amsopn}
\usepackage{amscd}

\usepackage [all]{xy}

\input xy

\newcommand{\sC}{\mathcal{C}}
\newcommand{\cinf}{\sC^\infty}

\newcommand{\sL}{\mathcal{L}}

\newcommand{\bR}{\mathbf{R}}

\newcommand{\bZ}{\mathbf{Z}}
\newcommand{\A}{\mathbb{A}}

\newcommand{\ruu}{\bR^{1|1}}
\newcommand{\rou}{\bR^{0|1}}

\newcommand{\tensor}{\otimes}

\newcommand{\tic}{\tilde{c}}

\newcommand{\nab}{\bar{\na}}
\newcommand{\Ab}{\bar{\A}}

\theoremstyle{plain}

\theoremstyle{definition}

\theoremstyle{remark}

\let\sym Q
\let\str D

\let\th\theta

\let\i\iota

\let\o\omega

\let\G\Gamma

\let\O\Omega

\let\na\nabla
\let\ra\rightarrow

\newcommand{\vd}{\partial_\th+\th\partial_t}

\begin{document}
\title{Addendum to \\
``Superconnections and Parallel Transport"}
\author{Florin Dumitrescu}
\date{\today}
\maketitle

\begin{abstract} In this addendum to our article ``Superconnections and Parallel Transport" we give an alternate  construction to the parallel transport of a superconnection contained in Corollary 4.4 of \cite{D1}, which has the advantage that is independent on the various ways  a superconnection splits as a connection plus a bundle endomorphism valued form.
\end{abstract}

\vspace{.2in}

Consider as in Section 4 of \cite{D1} a  superconnection $\A$ in the sense of Quillen (see \cite{Q} and \cite{BGV}) on a $\bZ/2$-graded vector bundle $E$ over a \emph{manifold} $M$, i.e. an odd first-order differential operator 
\[ \A: \O^*(M,E)\ra \O^*(M,E) \]
satisfying Leibniz rule
 \[ \A(\o\tensor s)= d\o\tensor s \pm \o\tensor \A(s), \]
with $\o\in\O^*(M)$ differential form on $M$ and $s\in \G(M;E)$ arbitrary section of the bundle $E$ over $M$. For such a superconnection we defined in \cite{D1} a notion of parallel transport along (families of) superpaths $c:S\times \ruu\ra M$ that is compatible under glueing of superpaths. Let us briefly recall this construction. First, let us write $\A=\A_1+A$, with $\A_1=\na$ the connection part of the superconnection $\A$ and $A\in\O^*(M, End\ E)^{odd}$ the linear part of the superconnection. For an arbitrary superpath $c$ in $M$ consider the diagram
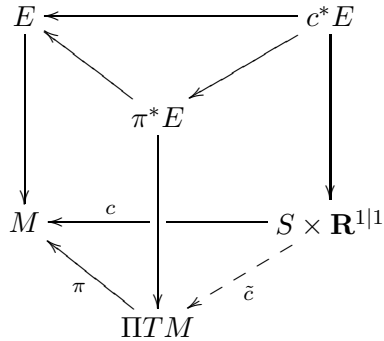
\[ \xymatrix{ E \ar[dd] & & c^*E \ar[ll] \ar[dl] \ar[dd] \\
& \pi^*E \ar[ul] \ar[dd] & \\
M & & S\times\ruu \ar'[l][ll]_{\ \ \ \ c} \ar@{-->}[dl]^-{\tilde{c}} \\
& \Pi TM \ar[ul]^\pi &  } \]
with $\tic$ a canonical lift of the path $c$ to $\Pi TM$, the ``odd tangent bundle" of $M$. Then parallel transport along $c$ is defined by {\it parallel} sections $\psi\in\G(c^*E)$ along $c$ which are solutions to the following differential equation 
\[ (c^*\na)_{D}\psi- (\tilde{c}^*A)\psi= 0.  \]
Here $D=\vd$ denotes the standard (right invariant) vector field on $\ruu$, see Section 2.4 of \cite{D1}.\\

Our alternate construction goes as follows. We first write $\A=\A_0+\Ab$, where $\A_0$ denotes the zero part of the superconnection and $\Ab$ the remaining part. Define then a connection $\nab$ on the bundle $\pi^*E$ over $\Pi TM$ as follows 
\[ \nab_{\sL_X}(\o\tensor s)= \sL_X\o\tensor s \pm \i_X\Ab s, \]
\[ \nab_{\i_X} (\o\tensor s)= \i_X\o\tensor s, \]
for $\o\in\O^*(M)$ and $s\in \G(M;E)$. Here, for a vector field $X$ on $M$, $\sL_X$ and $\i_X$ denote the Lie derivative respectively contraction in the $X$-direction acting as even respectively odd derivations on $\O^*(M)=\cinf(\Pi TM)$, i.e. as vector fields on $\Pi TM$. These relations are enough to define a connection $\nab$ on the bundle $\pi^*E$ over $\Pi TM$ since the algebra of vector fields on $\Pi TM$ is generated over $\cinf(\Pi TM)$ by vector fields of the type $\sL_X$ and $\i_X$, for $X$ arbitrary vector field on $M$, i.e.
\[ Vect(\Pi TM)= \cinf(\Pi TM)< \sL_X, \i_X\ | \ X\in Vect(M) > .\]
Parallel transport along a superpath $c:S\times \ruu\ra M$ is defined by {\it parallel sections}  $\psi\in\G(c^*E)$ along $c$ which are solutions to the following differential equation 
\[ (\tic^*\nab)_{D}\psi- (c^*\A_0)\psi= 0.  \]
As before, the parallel transport is well-defined (cf. Proposition 4.2 of \cite{D1}) by this ``half-order" differential equation and is compatible under glueing of superpaths (i.e. it satisfies properties (i) and (ii) of Theorem 4.3 in \cite{D1}). The advantage of this construction resides in the fact that the parallel transport so defined is invariant under the various ways in which a superconnection can be written as a sum of a connection plus a linear part, as $\Ab$ is invariant under such splittings. \\

\noindent Denote by $D$ the de Rham differential on $\Pi TM$. If $\o$ is a function on $\Pi TM$, then the 1-form $D\o$ on $\Pi TM$ evaluated on the standard odd vector field $d$ on $\Pi TM$ gives us
\[ (D\o)(d)= d\o, \]
the differential of $\o$, understood as a function on $\Pi TM$. This allows us to conclude that 
for any $s$ a section of $E$,
\[ \nab_d s= \Ab s. \]
Let us note that the connection $\nab$ is torsion free in the odd directions, i.e.
\[ [\nab_{\i_X},\nab_{\i_Y}]= \nab_{[\i_X,\i_Y]}\ (=0), \]
where $X$ and  $Y$ are vector fields on $M$. \\

\noindent The two constructions coincide when we consider connections instead of superconnections on the bundle $E$ over $M$. When the manifold $M$ reduces to a point, a graded vector bundle with superconnection reduces to a $\bZ/2$-vector space $V$ together with an odd endomorphism $A\ (=\A_0)$ of $V$. The two constructions of parallel transport we considered also coincide in this situation giving rise to the supergroup homomorphism of Example 3.2.9 in \cite{ST}
\[ \ruu \ni(t, \th)\longmapsto e^{-tA^2+\th A}\in GL(V), \]
as solution to the ``half-order" differential equation $D\psi= A\psi$.\\

\noindent We can also {\it recover} the superconnection from its associated parallel transport by first recovering the $\A_0$ from looking at constant superpaths in $M$, and then recover $\Ab$ by looking at parallel transport along the superpath given by $\ruu\times \Pi TM \ra \rou\times \Pi TM\ra M$, where the two maps are the obvious projection maps. The lift of such a superpath to $\Pi TM$ gives, after the obvious projection map, the flow of the vector field $d$ on $\Pi TM$. Given that $\nab_d s= \Ab s$, this recovers $\Ab$. Compare with Section 4.4 of \cite{D1}.

\vspace{.2in}

\noindent\emph{Acknowledgements.} The construction presented here is a mere continuation of an idea of Stephan Stolz who first thought to interpret a Quillen superconnection on a bundle $E$ over $M$ as a connection on the pullback bundle $\pi^*E$ over $\Pi TM$. I would like to thank Peter Teichner for suggesting to write up this Addendum.

\bibliographystyle{plain}
\bibliography{bibliografie}

\bigskip
\raggedright Institute of Mathematics of the Romanian Academy ``Simion Stoilow"\\  21 Calea Grivitei \\
Bucharest, Romania.\\ Email: {\tt florinndo@gmail.com}

\end{document}